\documentclass[centertags,a4, 11pt]{article}

\usepackage{amsmath,amsthm,amscd,amssymb}

\usepackage{tikz,subfigure,overpic}

\usetikzlibrary{shapes,snakes}

\usepackage{latexsym}

\usepackage{mathrsfs}
\usepackage{graphicx,subfig}
\usepackage{picins}
\usepackage{rotating}

\newcommand{\bbN}{{\mathbb{N}}}
\newcommand{\N}{{\mathbb{N}}}

\newcommand{\R}{{\mathbb{R}}}

\newcommand{\bbC}{{\mathbb{C}}}
\newcommand{\C}{{\mathbb{C}}}

\newcommand{\no}{\nonumber}

\newcommand{\supp}{\text{\rm{supp}}}

\newcommand{\beq}{\begin{equation}}

\newcommand{\eeq}{\end{equation}}

\newcommand{\ba}{\begin{align}}

\newcommand{\ea}{\end{align}}

\newcommand{\vphi}{\varphi}

\allowdisplaybreaks

\numberwithin{equation}{section}

\newtheorem{theorem}{Theorem}[section]

\newtheorem{proposition}[theorem]{Proposition}

\newtheorem{lemma}[theorem]{Lemma}

\newtheorem{corollary}[theorem]{Corollary}

\theoremstyle{definition}

\newtheorem{definition}[theorem]{Definition}

\newtheorem{example}[theorem]{Example}

\theoremstyle{remark}

\title{Spectral Analysis of Certain Spherically Homogeneous Graphs}
\author{Jonathan Breuer\footnote{ Einstein Institute of Mathematics, The Hebrew University of Jerusalem,
Jerusalem 91904, Israel. E-mail: jbreuer@math.huji.ac.il. Supported in part by The Israel Science Foundation (grant no.\ 1105/10)} \ and
Matthias Keller\footnote{Mathematisches Institut, Friedrich Schiller Universit\"at Jena,
  D-07743 Jena, Germany, Email: m.keller@uni-jena.de. Supported in part by The Israel Science Foundation (grant no.\ 1105/10)}}

\begin{document}
\date{}
\sloppy
\maketitle

\begin{abstract}
We study operators on rooted graphs with a certain spherical homogeneity. These graphs are called path commuting and allow for a decomposition of the adjacency matrix and the Laplacian into a direct sum of Jacobi matrices which reflect the structure of the graph. Thus, the spectral properties of the adjacency matrix and the Laplacian can be analyzed by means of the elaborated theory of Jacobi matrices.
For some examples which include antitrees, we derive the decomposition explicitly and present a zoo of spectral behavior induced by the geometry of the graph. In particular, these examples show that spectral types are not at all stable under rough isometries.
\end{abstract}

\section{Introduction}

Let $H$ be a bounded, selfadjoint operator acting on a Hilbert space $\mathcal{H}$. Let $\psi \in \mathcal{H}$, $\|\psi\|= 1$, and consider the restriction of $H$ to the cyclic subspace spanned by $H$ and $\psi$, $\mathcal{H}_\psi$. The set $\{\psi, H\psi, H^2 \psi, \ldots \}$ spans $\mathcal{H}_\psi$, so, by `Gram-Schmidting' it, one obtains an orthonormal basis for $\mathcal{H}_\psi$: $\{\psi, p_1(H)\psi, p_2(H)\psi, p_3(H)\psi, \ldots \}$ where $p_n(H)$ is a polynomial of degree $n$
in $H$. Clearly, $\langle p_j(H)\psi ,p_k(H)\psi\rangle=\delta_{j,k}$ (where $\langle\cdot, \cdot\rangle$ is the inner product in $\mathcal{H}$), and, since the subspace spanned by $\{\psi, H\psi, H^2\psi, \ldots, H^n \psi\}$ is equal to the subspace spanned by $\{\psi, p_1(H)\psi, p_2(H) \psi, \ldots, p_n(H) \psi \}$ for any $n$, it easily follows that the matrix of the restriction of $H$ to $\mathcal{H}_\psi$, in the basis  $\{\psi, p_1(H)\psi, p_2(H)\psi, p_3(H)\psi, \ldots \}$, is tridiagonal. Since, by Zorn's Lemma, $\mathcal{H}$ can be decomposed as a direct sum of spaces that are cyclic for $H$, see \cite{rs}, this implies that $H$ is equivalent to a direct sum of tridiagonal matrices. Such tridiagonal matrices are also called Jacobi matrices.

The (standard) argument sketched above implies that a natural strategy for the spectral analysis of selfadjoint operators is to decompose them as tridiagonal matrices, and then apply methods and results from the extensive theory of tridiagonal matrices. However, the process of tridiagonalization is usually a hard problem and only seldom can one obtain sufficient structural information about the resulting matrices.

One purpose of this paper is to describe a class of graphs such that the tridiagonalization and decomposition described above, for the associated adjacency matrix and Laplacian, are, in a certain sense, simple. As we shall see, this class includes some graphs that have been studied recently in other contexts. The second purpose of this paper is to then apply the strategy described above to obtain interesting spectral information about these graphs.

An example of such graphs is given by antitrees. These graphs serve as threshold examples to show disparities of phenomena in the discrete and continuous case. In particular, this concerns the relation of volume growth and properties of the heat flow on the one hand \cite{Hu,GHM,KLW, Woj} and the bottom of the (essential) spectrum on the other hand \cite{KLW}.

A simple version of the decomposition we shall describe below was applied in \cite{AF,Br,breuer-mol,BF,GG} to analyze spherically symmetric rooted trees. This work grew out of the realization that the tree structure is irrelevant for this strategy.
Moreover, as we shall show below, at least in the case of the adjacency matrix, even spherical symmetry is not the right type of symmetry to consider. Rather, we shall show that the right type of symmetry is a certain invariance of the number of paths of a given length between two vertices in the same generation, under changing the order of steps. We call graphs with this type of symmetry, \emph{path commuting}.

Path commuting graphs are described in Section~\ref{s:PCRG}. In this introduction we focus on a subclass of path commuting graphs that we call \emph{family preserving graphs}. Let $G$ be a rooted graph and let $S_n$ denote the sphere of radius $n$ in $G$, (that is, $S_n$ is the set of vertices of distance $n$ from the root).

\begin{definition} \label{f-preserving}
A rooted graph, $G$, is called \emph{family preserving} if the following three conditions hold:\\
(i) If $x, y \in S_n$ are both connected to a vertex $z \in S_{n+1}$, then there is a rooted graph automorphism, $\tau$, such that $\tau(x)=y$ and $\tau \upharpoonright S_{n+j}=Id$ for any $j \geq 1$. \\
(ii) If $x, y \in S_n$ are both connected to a vertex $z \in S_{n-1}$, then there is a rooted graph automorphism, $\tau$, such that $\tau(x)=y$ and $\tau \upharpoonright S_{n-j}=Id$ for any $1 \leq j \leq n$.\\
(iii) If $x, y \in S_n$ are neighbors (namely, there is an edge connecting $x$ and $y$), then there is a rooted graph automorphism, $\gamma$, such that $\gamma(x)=y$ and $\gamma(y)=x$.
\end{definition}

Figuratively speaking, a family preserving graph is a graph in which members of a single family `look alike' in the sense that their relationships to members of other `families' are also similar. Note that in the case of a rooted tree, family preserving is equivalent to spherical symmetry. Theorem~\ref{family-preserving} says that any family preserving graph is path commuting and Theorem~\ref{path-commuting} says that path commuting graphs admit a `nice' decomposition for the Laplacian and adjacency matrix, with the algorithm for this decomposition given in the proof.

The  machinery described here allows us to construct a zoo of examples of graphs with various types of spectral behavior resulting from the geometry of the graph.
In the following two subsections we present two examples. Both of these examples have been studied at various other places in the literature. One class is that of antitrees mentioned above and studied in \cite{Hu,GHM,KLW,Web,Woj}. One particular consequence of our method is the result that, in the context of antitrees, spectral types are not at all stable under rough isometries. Another class is that of spherically symmetric trees where occasionally edges are added such that spheres become complete graphs. Special cases of such models were studied for example in \cite{FHS,K}.

Let us next introduce the operators we are concerned with. To this end let $G$ be a graph that is locally finite, i.e.,   for each vertex $x$ the degree  $\deg(x)$, that is number of neighbors, is finite. We study the Laplace operator $\Delta$ on $\mathcal{H}\equiv\ell^{2}(G)$ given by
\begin{align*}
\Delta \phi(x)=\sum_{y\sim x}(\phi(x)-\phi(y)),
\end{align*}
where $x\sim y$ means that $x$ and $y$ are neighbors. Note that $\Delta$ is unbounded whenever there is no uniform bound on the vertex degree. However, $\Delta$ is essentially selfadjoint on the functions of compact (and thus finite) support and its domain is given by $D(\Delta)=\{\phi\in\ell^{2}(G)\mid \Delta \phi\in\ell^{2}(G)\}$, see \cite{KL,Woj1}. So, in the following $\Delta$ denotes the selfadjoint operator on $D(\Delta)$.

In contrast, the adjacency matrix $A$ given by
\begin{align*}
    A\phi(x)=\sum_{y\sim x}\phi(y),
\end{align*}
is not necessarily essentially selfadjoint, see \cite{M} and \cite{G} for recent developments. Therefore, for the sake of simplicity, we restrict ourselves to the case of bounded vertex degree whenever we consider the adjacency matrix. In this case, $A$ is a bounded symmetric operator and therefore selfadjoint.

For any vertex, $x \in G$, we let $\delta_x \in \mathcal{H}$ be the delta function at $x$. For a rooted graph, we let $s_n=\# S_n=$ the cardinality of $S_{n}$.

\subsection{Antitrees}\label{ss:antitrees}
\begin{figure}[!h]
\centering
\scalebox{0.5}{\includegraphics{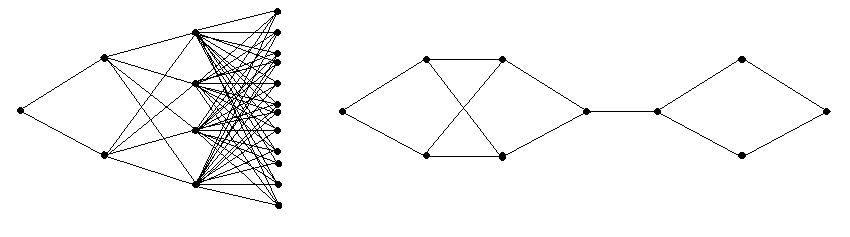}}
\caption{Two examples of  antitrees}\label{f:antitree}
\end{figure}

An \emph{antitree} is a connected rooted graph where every vertex in the $n$-th sphere, $n\geq1$, is connected to every vertex in the $(n-1)$-th and the $(n+1)$-th sphere and to none in the $n$-th sphere. See Figure~\ref{f:antitree} for two examples.

\begin{theorem}\label{t:FP_antitrees} Every antitree is family preserving.
\end{theorem}
\begin{proof} By definition, every vertex in a sphere is connected to all vertices in the previous and succeeding sphere. Therefore, for any two vertices in the same sphere there is a rooted graph automorphism interchanging these two vertices and leaving every other vertex invariant.
\end{proof}

It follows from Theorem~\ref{path-commuting} and~\ref{family-preserving} that the Laplacian on antitrees can be decomposed as a direct sum of Jacobi matrices. We describe the details in Section~\ref{s:Proof_Applications}. Here, we state a few theorems that are corollaries of the decomposition. The proofs will also be given in Section~\ref{s:Proof_Applications}.

\begin{theorem}\label{c:antitrees1} The spectrum of the Laplacian on an antitree is given by the spectrum of one infinite  Jacobi matrix and eigenvalues $(s_{n-1}+s_{n+1})$, $n\ge1$, with finitely supported eigenfunctions.
\end{theorem}

A sequence $v\in\R^{\N}$ is called \emph{eventually periodic} if  there is $N\in\bbN$ and $q\in\bbN$ such that $v_{n+q}=v_{n}$ for $n\geq N$. The following is an analogue of \cite[Theorem 1]{BF} for antitrees and like the theorem there, it is an immediate consequence of the decomposition and Remling's Theorem \cite{Rem}.

\begin{theorem}\label{c:antitrees2} Let an antitree  be given and let the cardinalities of the spheres $(s_{n})$ be bounded. Then, $\Delta$ has absolutely continuous spectrum if and only if $(s_{n})$ is eventually periodic.
\end{theorem}

The previous theorem shows that spectral types are  not at all stable under rough isometries. Recall that for two graphs $G_{1}$ and $G_{2}$  a rough isometry is a map $\Phi$ between the vertex sets such that there are constants $a,b,c\geq0$ such that $a^{-1} d_{1}(x,y)-b\leq d_{2}(\Phi(x),\Phi(y))\leq a d_{1}(x,y)+b$ for all $x,y\in G_{1}$ and such that for each $z \in G_2$ there is $x\in G_{1}$ with $d_{2}(\Phi(x),z)\leq c$, (where $d_{1}$ and $d_{2}$ are the natural graph metrics on $G_{1}$ and $G_{2}$.)
This instability is illustrated by the following example.

\begin{example}
Let an antitree $G$ be given which is determined by a bounded  sequence $(s_{n})$ that is not eventually periodic. Moreover, consider $\bbN_{0}$ as a graph with edges connecting $n$ and $n+1$ for all $n\geq1$. The two graphs are roughly isometric (i.e., one can choose $\Phi:G\to\bbN_{0}$ such every vertex in the $n$-th sphere is mapped to $n$, $n\geq0$).
However, it is well known that the Laplacian on $\bbN_{0}$ has purely absolutely continuous spectrum while by the theorem the Laplacian on $G$ has purely singular spectrum.
\end{example}

\subsection{Trees with complete spheres}\label{ss:treeswithcompletespheres}

Let $T=T(k)$ be a spherically  symmetric tree whose branching numbers are given by a sequence $k\in\N^{\N}$, i.e., the number of forward neighbors of a vertex in the sphere $S_{n}$ is $k_{n+1}$, $n\geq0$. For a sequence $\gamma\in\{0,1\}^{\N}$ let $G=G(k,\gamma)$ be the graph that results from $T(k)$ by connecting all vertices in a sphere $S_{n}$ by edges whenever $\gamma_{n}=1$, $n\geq 1$.

\begin{center}
\begin{figure}[!h]
\includegraphics[scale=0.6]{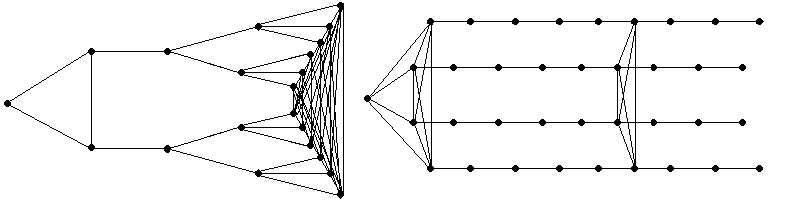}
\caption{Two examples of trees with complete spheres.}\label{f:tree}
\end{figure}
\end{center}

In \cite{FHS} the Anderson model is studied on a graph $G=G(k,\gamma)$, where $k\equiv 2$ and in every sphere edges with weight $\gamma\equiv \mathrm{const}\in(0,\infty)$ are inserted.  For the sake of simplicity, we deal here with unweighted edges only.

Moreover, the absence of essential spectrum of $\Delta$ on $G=G(k,1)$ with growing $k$ was studied in \cite{K}, (see also the application section of \cite{KLW}).

\begin{theorem}\label{t:FP_antitrees} For every  $k\in\N^{\N}$, $\gamma\in\{0,1\}^{\N}$, the graph $G(k,\gamma)$ is family preserving.
\end{theorem}
\begin{proof}Clearly $T(k)$ is family preserving. Moreover, a mapping on the vertices that is a rooted graph automorphism for $T(k)$ is a rooted graph automorphism for $G(k,\gamma)$ and vice versa. Thus, the statement follows.
\end{proof}

Again, the decomposition into the relevant Jacobi matrices is described in in Section~\ref{s:Proof_Applications}. Here, we state two theorems that are corollaries and are also proven there. We first present some examples of  graphs with a bounded absolutely continuous spectral component and infinitely many discrete eigenvalues.

\begin{theorem}\label{c:treeswithcompletespheres1}Let  $k$ be eventually periodic and $\gamma\equiv 1$. Then, the spectrum of $\Delta$ on $G(k,\gamma)$ consists of finitely many bands of absolutely continuous spectrum of multiplicity one and discrete eigenvalues accumulating at infinity.
\end{theorem}

The next theorem gives an example of a graph with one absolutely continuous and finitely many singular continuous components.
\begin{theorem}\label{c:treeswithcompletespheres2}
Let $\kappa\geq2$ and  $L_{j}\geq \prod_{i=1}^{j}(\kappa-1+\log i)$, $j\geq1$. Let  $k_{1}=\kappa$, $k_{n}=1$, $n\geq2$, and $\gamma_{n}=1$ if $n=\sum_{j=1}^{k}L_{j}$ for some $k\geq1$ and $\gamma_{n}=0$ otherwise. Then the spectrum of $\Delta$ on $G(k,\gamma)$ has an absolutely continuous component of multiplicity one and a singular continuous component of multiplicity $(\kappa-1)$, both of which consist of the interval $[0,4]$. The point $4$ may be an eigenvalue of finite multiplicity and an accumulation point for eigenvalues outside $[0,4]$. Moreover, the point $(2+\sqrt{4+\kappa^2})$ is a point of accumulation of eigenvalues (and perhaps an eigenvalue as well).
\end{theorem}

\bigskip

We note that while we restrict our attention in this paper to \emph{rooted} graphs, our methods seem to extend to unrooted graphs as well. The essential requirement is that of a direction on a graph. Being rooted simply means that the graph ends in one direction. Thus, it seems natural to also consider decompositions into two sided Jacobi matrices to treat `two sided' path commuting graphs. For the sake of simplicity, we will not do this here. However, it seems that some of the admissible graphs of \cite{MRT} (though probably not all) are path commuting and may be treated by the methods described in our paper.

The rest of this paper is structured as follows. Path commuting rooted graphs are described in Section~\ref{s:PCRG} where the relevant  decomposition is also given as Theorem~\ref{path-commuting}. As the characterization of path commuting graphs is unfortunately not local, we want to consider graphs that are simpler to describe. These are the family preserving graphs defined above. Theorem~ \ref{family-preserving} states that family preserving graphs are path commuting. The proofs of the statements of Section~\ref{s:PCRG} are given in Section~\ref{s:Proofs} and the proofs of the results described in Subsections~\ref{ss:antitrees} and \ref{ss:treeswithcompletespheres} are given in Section~\ref{s:Proof_Applications}.

\section{Path Commuting Rooted Graphs}\label{s:PCRG}

Throughout this paper we shall be dealing with rooted, locally finite, connected, simple graphs. A rooted graph is a  graph, $G$, with a distinguished vertex, $o$, that we call the root. Simplicity means there are no multiple edges and no self loops. The existence of the root induces a natural ordering on the vertices according to spheres, that is, with respect to their distance to their root. We let $S_n$ denote the sphere of distance $n$ to the root, namely the set of vertices $x$ such that $d(x,o)=n$. We shall sometimes call $S_n$ the $n$'th generation. We do not distinguish in notation between $G$ and its vertex set.

A rooted graph automorphism is an automorphism that fixes the root. It follows that if $\tau$ is a rooted graph automorphism of $G$ and $x \in S_n \subseteq G$, then $\tau(x) \in S_n$ as well. A path of length $n$ in $G$ is an $n$-tuple of vertices $(x_0, x_1, \ldots, x_n )$ such that there is an edge connecting $x_j$ to $x_{j+1}$ for all $1 \leq j  \leq n$. We need some simple definitions.

\begin{definition}
Let $G$ be a rooted graph and let $x,y \in S_n$.\\
A \emph{$k$-forward path} from $x$ to $y$ is a path $(x_0,x_2, \ldots, x_{2k})$  of length $2k$ such that $x_0=x$, $x_{2k}=y$ and such that $x_j \in S_{n+j}$ for all $j \leq k$ (which implies that $x_j \in S_{n+2k-j}$ for $j \geq k$).

Similarly, a \emph{$k$-backward path} from $x$ to $y$ is a path $(x_0,x_2, \ldots, x_{2k})$ of length $2k$ such that $x_0=x$, $x_{2k}=y$ and such that $x_j \in S_{n-j}$ for all $j \leq k$ (which implies that $x_j \in S_{n-2k+j}$ for $j \geq k$).
\end{definition}

In simple words, a $k$-forward path from $x$ to $y$ is a path that takes $k$ steps away from the root (``forward'') and then $k$ steps towards the root (``backward'') in order to go from $x$ to $y$. A backward path is the same with the backward steps taken first. Note that forward or backward paths starting at $x\in S_n$ always end up at some vertex in $S_n$, by definition.

\begin{definition}
Let $G$ be a rooted graph and let $x,y \in S_n$.

A \emph{$k, \ell$-forward-backward path} ($(k, \ell)$-f.b.\ path, for short) from $x$ to $y$ is a $k$-forward path, starting at $x$ and ending at some vertex $z \in S_n$, followed by an $\ell$-backward path starting at $z$ and ending at $y$.

Similarly, a \emph{$k, \ell$-backward-forward path} ($(k,\ell)$-b.f.\ path, for short) from $x$ to $y$ is a $k$-backward path, starting at $x$ and ending at some vertex $z \in S_n$, followed by an $\ell$-forward path starting at $z$ and ending at $y$.
\end{definition}
We also need to allow for steps to be taken within $S_n$.
\begin{definition}
Let $G$ be a rooted graph and let $x,y \in S_n$.\\
A \emph{tailed-$k$-forward path} (or tailed-$k$-f.\ path for short) from $x$ to $y$ is a path $(x_0,x_2, \ldots, x_{2k+1})$ of length $2k+1$, where $(x_0,x_2, \ldots, x_{2k})$ is a $k$-forward path from $x$ to $x_{2k}\in S_{n}$ and $x_{2k}\sim x_{2k+1}=y$. 

A  \emph{tailed-$k$-backward path} (or tailed-$k$-b.\ path for short) from $x$ to $y$ is a path $(x_0,x_2, \ldots, x_{2k+1})$ of length $2k+1$, where $(x_0,x_2, \ldots, x_{2k})$ is a $k$-backward path from $x$ to $x_{2k}\in S_{n}$ and $x_{2k}\sim x_{2k+1}=y$.

A \emph{headed-$k$-forward path} (or headed-$k$-f.\ path for short) from $x$ to $y$ is a path $(x_0,x_2, \ldots, x_{2k+1})$ of length $2k+1$, where $x\sim x_{1}\in S_{n}$ and $(x_1,x_2, \ldots, x_{2k+1})$ is a $k$-forward path from $x_{1}$ to $x_{2k+1}=y$.

A \emph{headed-$k$-backward path} (or headed-$k$-b.\ path for short) from $x$ to $y$ is a path $8x_0,x_2, \ldots, x_{2k+1})$  of length $2k+1$, where $x\sim x_{1}\in S_{n}$ and $(x_1,x_2, \ldots, x_{2k+1})$ is a $k$-backward path from $x_{1}$ to $x_{2k+1}=y$.
\end{definition}

In other words, \emph{tailed paths} are paths that end with a step within $S_n$ and \emph{headed paths} are paths that begin with a step within $S_n$.

Let $G$ be a rooted graph. For two vertices, $x, y \in S_n$, let
\begin{align*}
\#fb_{(k, \ell)}(x,y)&=\#\{\mbox{$(k,\ell)$-f.b.\ paths from $x$ to $y$}\},\\
\#bf_{(k,\ell)}(x,y)&=\#\{\mbox{$(k,\ell)$-b.f.\ paths from $x$ to $y$}\},\\
\widetilde{\#f}_k(x,y)&=\#\{\mbox{tailed-$k$-f.\ paths from $x$ to $y$}\},\\
\widetilde{\#b}_k(x,y)&=\#\{\mbox{tailed-$k$-b.\ paths from $x$ to $y$}\},\\
\widehat{\#f}_k(x,y)&=\#\{\mbox{headed-$k$-f.\ paths from $x$ to $y$}\},\\
\widehat{\#b}_k(x,y)&=\#\{\mbox{headed-$k$-b.\ paths  from $x$ to $y$}\}.
\end{align*}


\begin{definition}
A rooted graph, $G$, is called \emph{path commuting} if for any $n$, $x,y \in S_n$ and $k, \ell \geq 0$
\begin{align*}
\#fb_{(k, \ell)}(x,y)&=\#bf_{(\ell,k)}(x,y),\\
\widetilde{\#f}_k(x,y)&=\widehat{\#f}_k(x,y),\\
\widetilde{\#b}_k(x,y)&=\widehat{\#b}_k(x,y).
\end{align*}

A rooted graph is called \emph{strongly path commuting} if it is path commuting and in addition, for any $x, y \in S_n$,  $\deg(x)=\deg(y)$. Namely, the degree is a function of the distance from the root.
\end{definition}

To give a simple example of a (strongly) path commuting rooted graph, consider a rooted spherically symmetric tree, $\Gamma$. By \emph{spherically symmetric} we mean that for any two vertices on $S_j$, there is an automorphism sending one to the other. For rooted trees, this is equivalent to $x$ and $y$ having the same degree (for any fixed $j$). That $\Gamma$ is strongly path commuting is easy to check since $0=\widetilde{\#f}_k(x,y)=\widehat{\#f}_k(x,y) =\widetilde{\#b}_k(x,y)=\widehat{\#b}_k(x,y)$ for any $x,y \in S_n$ (since, by the tree property, there are no edges between vertices on the same sphere). Moreover, since there are no forward paths between two different vertices in $S_j$ (only from a vertex to itself), we obtain that $\#fb_{(k, \ell)}(x,y)$ is a product of the number of $k$-forward paths from $x$ to itself, and the number of $\ell$-backward paths from $x$ to $y$. (Note that the number of $\ell$-backward paths from $x$ to $y$ is always $0$ or $1$.) By the same reasoning, $\#bf_{(\ell,k)}(x,y)$ is a product of the number of $\ell$-backward paths from $x$ to $y$, and the number of $k$-forward paths from $y$ to itself. By the spherical symmetry, these are the same. Therefore, spherically symmetric rooted trees are strongly path commuting. In fact, the reverse is true as well. Clearly, strongly path commuting rooted trees are spherically symmetric. Moreover, even path commuting rooted trees are necessarily spherically symmetric,  since if $\Gamma$ is a rooted tree that is path commuting, then by considering only $(1,\ell)$-f.b.\ and $(\ell,1)$-b.f.\ paths, together with the fact that for any $x, y \in S_n$ there is a backward path going from $x$ to $y$ (if needed, one can always return all the way back to the root), we see that the number of forward neighbors is constant across spheres. Since the number of backward neighbors is $1$ for any vertex, we get the claim. Thus,

\begin{proposition} \label{trees}
A rooted tree, $\Gamma$, is path commuting iff it is strongly path commuting iff it is spherically symmetric.
\end{proposition}

The reason we single out (strongly) path commuting graphs, is that on these graphs there is natural way for decomposing $\Delta$ and $A$. In particular, $A$ decomposes on path commuting graphs, and $\Delta$ decomposes on strongly path commuting graphs. The extra assumption on the degrees is needed for $\Delta$ because of the additional diagonal terms it has.

\begin{theorem} \label{path-commuting}
(1.) Let $G$ be a path commuting graph with bounded degree. Then $\ell^2(G)=\oplus_{r=0}^\infty \mathcal{H}_r$ such that:
\begin{itemize}
  \item [(i)] For each $r$ there exists $n(r)$ and a vector $\phi_{0}^r$ such that $\supp(\phi_0^r) \subseteq S_{n(r)}$ and $\mathcal{H}_r =\overline{\textrm{span}\{\phi_0^r,A\phi_0^r, A^2 \phi_0^r, \ldots\}}$, that is, $\mathcal{H}_r$ is the cyclic subspace spanned by $\phi_0^r$ and $A$. In particular, $\mathcal{H}_r$ is $A$-invariant and $A$ is unitarily equivalent to the direct sum of its restrictions to $\mathcal{H}_r$.
  \item [(ii)]
 The set $\{\phi_0^r, \phi_1^r, \phi_2^r, \dots\}$, obtained from $\{\phi_0^r,A\phi_0^r, A^2 \phi_0^r, \ldots\}$ by applying the Gram-Schmidt process has the property that $\supp (\phi_k^r) \subseteq S_{n(r)+k}$ for each $k \geq 0$.
\end{itemize}
(2.) Let $G$ be a strongly path commuting graph. Then $\ell^2(G)=\oplus_{r=0}^\infty \mathcal{H}_r$ such that:
\begin{itemize}
\item [(i)] For each $r$ there exists $n(r)$ and a vector $\phi_{0}^r$ such that $\supp(\phi_0^r) \subseteq S_{n(r)}$ and $\mathcal{H}_r =\overline{\textrm{span}\{\phi_0^r,\Delta\phi_0^r, \Delta^2 \phi_0^r, \ldots\}}$, that is, $\mathcal{H}_r$ is the cyclic subspace spanned by $\phi_0^r$ and $\Delta$. In particular, $\mathcal{H}_r$ is $\Delta$-invariant and $\Delta$ is unitarily equivalent to the direct sum of its restrictions to $\mathcal{H}_r$.
  \item [(ii)] The set $\{\phi_0^r, \phi_1^r, \phi_2^r, \dots\}$, obtained from $\{\phi_0^r, \Delta \phi_0^r, \Delta^2 \phi_0^r, \ldots\}$ by applying the Gram-Schmidt process has the property that $\supp( \phi_k^r )\subseteq S_{n(r)+k}$ for each $k \geq 0$.
\end{itemize}
\end{theorem}

The proof of Theorem \ref{path-commuting} will be given in Section~3. In the meantime, to illustrate why this notion is useful, consider the case of trees. Consider the simple trees in Figure~3. The tree in Figure~3(a) is not spherically symmetric. It is easy to see that for that example, the cyclic subspace spanned by $A$ and $\delta_o$ includes the both the function $\delta_{v_2}+\delta_{v_3}$ and the function $\delta_{v_2}+2\delta{v_3}$. On the other hand, for the tree in Figure~3(b), it is easy to see that the cyclic subspace spanned by $A$ and $\delta_o$ includes only the spherically symmetric functions on spheres and thus the resulting Jacobi matrix has a straightforward relation to the structure of the tree.

\setlength{\unitlength}{1cm}
\begin{figure}[!h] \label{figure3}
\begin{center}
\subfigure[A non path commuting rooted tree]{
\begin{tikzpicture}
\put(-3,0){\line(1,0){1}}
\draw (-3,-0.3) node{$o$};
\put(-3,0){\circle*{0.15}}
\put(-2,0){\circle*{0.15}}
\draw(-2.1,-0.3) node{$v_1$};
\put(-2,0){\line(1,2){0.5}} 
\put(-1.5,1){\circle*{0.15}}
\draw(-1.8,0.9) node{$v_2$};
\put(-2,0){\line(1,-2){0.5}} 
\put(-1.5,-1){\circle*{0.15}}
\draw(-1.5,-1.3) node{$v_3$};
\put(-1.5,1){\line(1,0){1}} 
\put(-0.5,1){\circle*{0.15}}
\draw(-0.3,0.7) node{$v_4$};
\put(-0.5,-0){\circle*{0.15}}
\draw(-0.3,-0.3) node{$v_5$};
\put(-0.5,-1){\circle*{0.15}}
\draw(-0.3,-1.3) node{$v_6$};
\put(-1.5,-1){\line(1,1){1}}
\put(-1.5,-1){\line(1,0){1}} 
\end{tikzpicture}}
\hspace*{3cm}
\subfigure[A path commuting rooted tree]{
\begin{tikzpicture}
\put(1,0){\line(1,0){1}}
\put(1,0){\circle*{0.15}}
\draw(1,-0.3) node{$o$};
\put(2,0){\circle*{0.15}}
\draw(1.9,-0.3) node{$v_1$};
\put(2,0){\line(1,2){0.5}}
\put(2.5,1){\circle*{0.15}}
\draw(2.2,1) node{$v_2$};
\put(2,0){\line(1,-2){0.5}}
\put(2.5,-1){\circle*{0.15}}
\draw(2.5,-1.3) node{$v_3$};
\put(2.5,1){\line(1,0){1}}
\put(2.5, -1){\line(1,0){1}}
\put(3.5,1){\circle*{0.15}}
\draw(3.5,0.7) node{$v_4$};
\put(3.5, -1){\circle*{0.15}}
\draw(3.5,-1.3) node{$v_5$};
\end{tikzpicture}}
\end{center}
\caption{Two illustrative examples}
\end{figure}
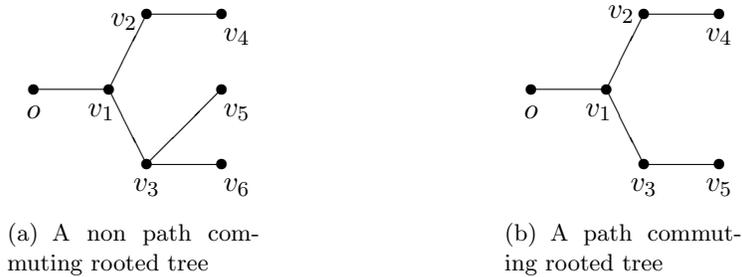
It is intuitively obvious that a path commuting graph has a certain amount of symmetry. One might think that spherical symmetry is either necessary or sufficient. Figure~4 shows neither is true.
The graph in Figure~4(a) is clearly not spherically symmetric ($v_4$ and $v_5$ are both in $S_3$ but there is no graph automorphism taking $v_4$ to $v_5$) but it is easy to check that it is path commuting. The graph in Figure~4(b), on the other hand, is spherically symmetric. However, it is easy to see that it is not path commuting: In $S_3$, $v_7$ and $v_5$ are connected by a $(1,1)$-b.f.\ path but there is no $(1,1)$-f.b.\ path between them. We note that the graph in Figure~4(a) is \emph{not} strongly path commuting and so also serves as an example to show that not all path commuting graphs are strongly path commuting.
\setlength{\unitlength}{1cm}
\begin{figure}[!h]
\begin{center}
\subfigure[Path commuting but not spherically symmetric]{
\begin{tikzpicture}
\put(-3,0){\line(1,0){1}}
\draw (-3,-0.3) node{$o$};
\put(-3,0){\circle*{0.15}}
\put(-2,0){\circle*{0.15}}
\draw(-2.1,-0.3) node{$v_1$};
\put(-2,0){\line(1,1){1}} 
\put(-1,1){\circle*{0.15}}
\draw(-1.3,1.3) node{$v_2$};
\put(-2,0){\line(1,-1){1}} 
\put(-1,-1){\circle*{0.15}}
\draw(-1.3,-1.3) node{$v_3$};
\put(-1,1){\line(1,0){1}}
\put(0,1){\circle*{0.15}}
\draw(-0,1.3) node{$v_4$};
\put(-0,0){\circle*{0.15}}
\put(-1,1){\line(1,-1){1}}
\draw(-0,-0.4) node{$v_5$};
\put(-0,-1){\circle*{0.15}}
\draw(-0,-1.3) node{$v_6$};
\put(-1,-1){\line(1,1){1}}
\put(-1,-1){\line(1,0){1}} 
\put(0.3,1){$\ldots$}
\put(0.3, -0){$\ldots$}
\put(0.3, -1){$\ldots$}
\put(1, -1){ }
\end{tikzpicture}}
\hspace*{1cm}
\subfigure[Spherically symmetric but  not path commuting]{
\begin{tikzpicture}
\put(1,1.5){\line(1,0){1}} 
\put(1,1.5){\circle*{0.15}} 
\draw(.7,1.4) node{$v_1$};
\put(2,1.5){\circle*{0.15}} 
\draw(1.9,1.2) node{$o$};
\put(2,1.5){\line(1,2){0.5}}
\put(2.5,2.5){\circle*{0.15}}
\draw(2.2,2.5) node{$v_2$};
\put(2,1.5){\line(1,-2){0.5}} 
\put(2.5,0.5){\circle*{0.15}}
\draw(2.2,0.5) node{$v_3$};
\put(2.5,0.5){\line(-1,-2){.5}}
\put(2.5,2.5){\line(1,0){1}} 
\put(2.5, 0.5){\line(1,0){1}} 
\put(3.5,2.5){\circle*{0.15}}
\draw(3.9,2.4) node{$v_7$};
\put(3.5,2.5){\line(1,-2){0.5}} 
\put(4,1.5){\circle*{0.15}} 
\put(3.5, 0.5){\circle*{0.15}} 
\draw(3.9,0.5) node{$v_8$};
\put(3.5,0.5){\line(1,2){0.5}} 
\draw(4.7,1.5)node{$v_{12}\ldots$};
\put(2, -.5){\circle*{0.15}} 
\draw(2.5, -0.5)node{$v_9$};
\put(2,-0,5){\line(-1,0){1}}
\put(2.5,2.5){\line(-1,2){.5}} 
\put(2,3.5){\circle*{0.15}} 
\draw(2.5,3.5)node{$v_6$};
\put(2,3.5){\line(-1,0){1}} 
\put(1,1.5){\line(-1,2){0.5}} 
\put(0.5, 2.5){\circle*{0.15}}
\draw(0.2, 2.5)node{$v_5$};
\put(0.5,2.5){\line(1,2){0.5}}
\put(1,3.5){\circle*{0.15}}
\draw(.5, 3.5)node{$v_{11}$};
\draw(.5, 4)node{$\ddots$};
\put(1,1.5){\line(-1,-2){0.5}}
\put(0.5,0.5){\circle*{0.15}}
\draw(0.2,0.5)node{$v_4$};
\put(1,-0.5){\line(-1,2){0.5}}
\put(1, -0.5){\circle*{0.15}}
\draw(0.5, -0.5)node{$v_{10}$};
\draw(.5,-0.8)node{\text{\reflectbox{$\ddots$}}};
\end{tikzpicture}}
\end{center}
\caption{Path commuting versus spherically symmetric. The dots signify some simple continuation (say, by copies of $\mathbb{N}$).}
\end{figure}
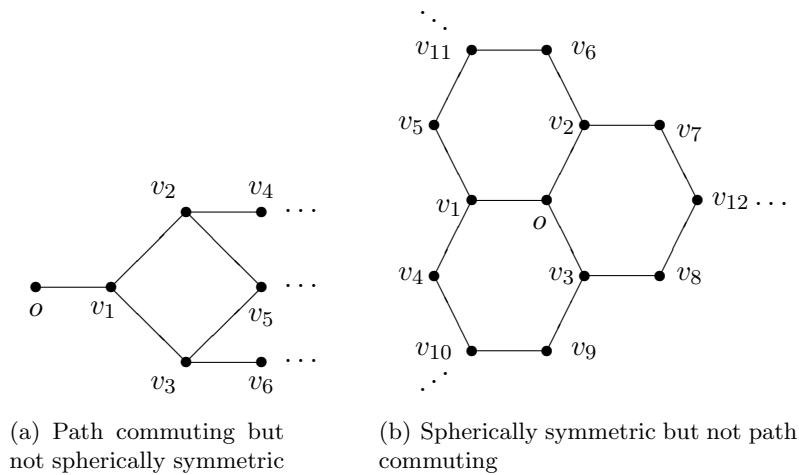

The notion of path commuting graph is a difficult notion to check since, in principle, it requires counting an infinite number of paths. We therefore restrict our attention to a different, \emph{local}, notion which is easier to check. That is the notion of family preserving graphs (recall Definition \ref{f-preserving}). We repeat the definition below.

\begin{definition}
We say that two vertices, $x, y \in S_n$ of a rooted graph $G$, are \emph{forward brothers}, if they have a common neighbor in $S_{n+1}$.

We say that two vertices, $x, y \in S_n$ are \emph{backward brothers}, if they have a common neighbor in $S_{n-1}$.
\end{definition}

\begin{definition}
A rooted graph, $G$, is called \emph{family preserving} if the following three conditions hold:\\
(i) If $x, y \in S_n$ are forward brothers, then there is a rooted graph automorphism, $\tau$, such that $\tau(x)=y$ and $\tau \upharpoonright S_{n+j}=Id$ for any $j \geq 1$.  \\
(ii) If $x, y \in S_n$ are backward brothers, then there is a rooted graph automorphism, $\tau$, such that $\tau(x)=y$ and $\tau \upharpoonright S_{n-j}=Id$ for any $1 \leq j \leq n$. \\
(iii) If $x, y \in S_n$ are neighbors, then there is a rooted graph automorphism, $\gamma$, such that $\gamma(x)=y$ and $\gamma(y)=x$.
\end{definition}

The following theorem will also be proven in the next section.
\begin{theorem}\label{family-preserving}
Any family preserving graph is strongly path commuting.
\end{theorem}

We shall show that any family preserving graph is spherically symmetric (see Lemma~\ref{spherical-symmetry}). Thus, it follows from the example in Figure 2(a) that there are path commuting graphs that are not family preserving. There are also strongly path commuting graphs that are not family preserving. An example is given in Figure~5.
\setlength{\unitlength}{1cm}
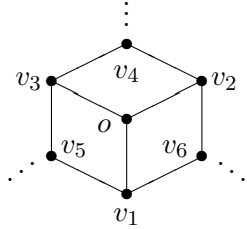
\begin{figure}[!h]\label{figure5}
\begin{center}
{
\begin{tikzpicture}
\put(0,1.5){\circle*{0.15}}
\draw(-0.3,1.4) node{$o$};
\put(0,1.5){\line(0,-1){1}}
\put(0,.5){\circle*{0.15}}
\draw(0,.2) node{$v_1$};
\put(0,1.5){\line(2,1){1}}
\put(1,2){\circle*{0.15}}
\draw(1.3,2) node{$v_2$};
\put(0,1.5){\line(-2,1){1}}
\put(-1,2){\circle*{0.15}}
\draw(-1.3,2) node{$v_3$};
\put(1,2){\line(-2,1){1}}
\put(-1,2){\line(2,1){1}}
\put(0,2.5){\circle*{0.15}}
\draw(0,2.1) node{$v_4$};
\draw(0,3) node{$\vdots$};
\put(1,2){\line(0,-1){1}}
\put(1,1){\circle*{0.15}}
\draw(.65,1.1) node{$v_6$};
\draw(1.4,.9) node{$\ddots$};
\put(-1,2){\line(0,-1){1}}
\put(-1,1){\circle*{0.15}}
\draw(-.7,1.1) node{$v_5$};
\draw(-1.4,.9) node{\text{\reflectbox{$\ddots$}}};
\put(-1,1){\line(2,-1){1}}
\put(1,1){\line(-2,-1){1}}
\end{tikzpicture}}
\end{center}
\caption{A strongly path commuting graph that is not family preserving.}
\end{figure}

\section{Proof of Theorems \ref{path-commuting} and \ref{family-preserving}}\label{s:Proofs}

We start by  recalling the representation of the graph adjacency matrix in spherical coordinates \cite{FHS1}:
\beq \label{spherical-coordinates}
A=\left(
\begin{array}{ccccc}
V_0    & E_0^\top & 0      & \ldots      & \ldots \\
E_0    & V_1 & E_1^\top    & \ddots     & 0 \\
0      &  E_1 & V_2    & E_2^\top    & \vdots \\
\vdots     & \ddots & E_2    &\ddots     &\ddots \\
\vdots & \ddots   & \ddots & \ddots & \ddots \\
\end{array} \right).
\eeq

Here, $V_j$ is the adjacency matrix of $S_j$ (so $V_0=0$ since there are no loops) and the matrix $E_j: \ell^2(S_j) \rightarrow \ell^2(S_{j+1})$ has a $1$ in the $(x, y)$-position if there is an edge between $x \in S_{j+1}$ and $y \in S_{j}$. Otherwise, it has a $0$ there.

The representation for $\Delta$ is similar:
\beq \label{spherical-coordinates}
\Delta=\left(
\begin{array}{ccccc}
\widetilde{V}_0    & \widetilde{E}_0^\top & 0      & \ldots      & \ldots \\
\widetilde{E}_0    & \widetilde{V}_1 & \widetilde{E}_1^\top    & \ddots     & 0 \\
0      &  \widetilde{E}_1 & \widetilde{V}_2    & \widetilde{E}_2^\top    & \vdots \\
\vdots     & \ddots & \widetilde{E}_2    &\ddots     &\ddots \\
\vdots & \ddots   & \ddots & \ddots & \ddots \\
\end{array} \right),
\eeq
where $\widetilde{E}_j=-E_j$ and $\widetilde{V}_j=D_j-V_j$ where $D_j$ is a diagonal matrix whose  entry at $(x,x)$ is the number of neighbors of $x$.

We let
\begin{align*}
\Lambda_{n,+j}&=E_{n}^\top E_{n+1}^\top \ldots  E_{n+j}^\top E_{n+j} \ldots E_{n+1} E_n,\quad\quad n, j \geq 0,\\
\Lambda_{n,-j}&=E_{n-1}E_{n-2} \ldots E_{n-j} E_{n-j}^\top \ldots E_{n-2}^{\top}E_{n-1}^\top,\quad n \geq 1,j \leq n.
\end{align*}
Then $\Lambda_{n,\pm j}: \ell^2(S_n) \rightarrow \ell^2(S_n)$. For notational convenience, when $n=0$ or for $j > n$, we let $\Lambda_{n,-j}=Id : \ell^2(S_n) \rightarrow \ell^2(S_n)$. We use $\widetilde{\Lambda}$ for the corresponding matrices in the case of the Laplacian.

\begin{lemma} \label{commute}
A rooted graph $G$ is path commuting iff for any $n$, the matrices in the set $\{V_n, \Lambda_{n+j}, \Lambda_{n,-j} \mid j=0,1,\ldots \}$ all commute with each other. If $G$ is strongly path commuting, then the matrices in the set $\{\widetilde{V}_n, \widetilde{\Lambda}_{n+j}, \widetilde{\Lambda}_{n,-j} \mid j=0,1,\ldots \}$ all commute with each other.
\end{lemma}

\begin{proof}
Fix $n$ and consider $x, y \in S_n$. It is easy to see that for any $k , \ell$ with $\ell \leq n$,
\beq \no
\# fb_{k,\ell}(x,y)=\langle \delta_y, \Lambda_{n,-\ell}\Lambda_{n,+k} \delta_x \rangle.
\eeq
Similarly,
\beq \no
\# bf_{\ell, k}(x,y)=\langle \delta_y, \Lambda_{n,k}\Lambda_{n,-\ell} \delta_x \rangle,
\eeq
\beq \no
\widehat{\# f}_k(x,y)=\langle\delta_y, \Lambda_{n,+k}V_n \delta_x \rangle,
\eeq
\beq \no
\widetilde{\# f}_k(x,y)=\langle\delta_y, V_n \Lambda_{n,+k} \delta_x \rangle,
\eeq
\beq \no
\widehat{\# b}_k(x,y)=\langle\delta_y, \Lambda_{n,-k}V_n \delta_x \rangle
\eeq
and
\beq \no
\widetilde{\# b}_k(x,y)=\langle\delta_y, V_n \Lambda_{n,-k} \delta_x \rangle.
\eeq
Thus, the statement follows directly from the definition (the fact that $\Lambda_{n,+j}$ and $\Lambda_{n,+k}$ commute, as well as the fact that $\Lambda_{n,-j}$ and $\Lambda_{n,-k}$ commute, follows easily from the fact that the other operators commute).

Recalling that for strongly path commuting graphs, $D_j$ are scalar matrices, the same proof goes to show that the corresponding matrices commute.
\end{proof}

\begin{lemma} \label{commute2}
With the notation as in Lemma \ref{commute}, let $G$ be a rooted path commuting graph. Then for any $n, j\geq0$,  the matrix $\Lambda_{n,+j}$ commutes with $E_{n}^\top E_{n+1}^\top \ldots  E_{n+j}^\top V_{n+j} E_{n+j} \ldots E_{n+1} E_n$, and for any $n \geq 1$ and $j \leq n$, $\Lambda_{n,-j}$ commutes with $ E_{n-1} \ldots E_{n-j} V_{n-j} E_{n-j}^\top \ldots E_{n-1}^\top$.

If $G$ is strongly path commuting, the analogous statement holds with $\Lambda$, $E$ and $V$ replaced by $\widetilde{\Lambda}$, $\widetilde{E}$ and $\widetilde{V}$ respectively.
\end{lemma}

\begin{proof}
Write
\beq \no
\begin{split}
& (E_{n}^\top E_{n+1}^\top \ldots  E_{n+j}^\top V_{n+j} E_{n+j} \ldots E_{n+1} E_n) \Lambda_{n,+j} \\
&=  (E_{n}^\top E_{n+1}^\top \ldots  E_{n+j}^\top V_{n+j} E_{n+j} \ldots E_{n+1} E_n)(E_{n}^\top E_{n+1}^\top \ldots  E_{n+j}^\top E_{n+j} \ldots E_{n+1} E_n)\\
&= E_{n}^\top E_{n+1}^\top \ldots  E_{n+j}^\top V_{n+j} (E_{n+j} \ldots E_{n+1} E_n E_{n}^\top E_{n+1}^\top \ldots  E_{n+j}^\top) E_{n+j} \ldots E_{n+1} E_n \\
&= E_{n}^\top E_{n+1}^\top \ldots  E_{n+j}^\top (E_{n+j} \ldots E_{n+1} E_n E_{n}^\top E_{n+1}^\top \ldots  E_{n+j}^\top)V_{n+j}
E_{n+j} \ldots E_{n+1} E_n \\
&=(E_{n}^\top E_{n+1}^\top \ldots  E_{n+j}^\top E_{n+j} \ldots E_{n+1} E_n) (E_{n}^\top E_{n+1}^\top \ldots  E_{n+j}^\top V_{n+j}
E_{n+j} \ldots E_{n+1} E_n) \\
&= \Lambda_{n, +j} (E_{n}^\top E_{n+1}^\top \ldots  E_{n+j}^\top V_{n+j}
E_{n+j} \ldots E_{n+1} E_n),
\end{split}
\eeq
where the second equality holds by Lemma \ref{commute}. The proof of the second statement is similar.

The proof of the statement for strongly path commuting graphs is obtained from this proof by placing a tilde "$\sim$" over the appropriate operators.
\end{proof}

\begin{proof}[Proof of Theorem \ref{path-commuting}]
We start by proving the statement regarding the adjacency matrix on a path commuting graph with bounded degree.
Let $P_n: \ell^2(G) \rightarrow \ell^2(G)$ be the orthogonal projection with range $\ell^2(S_n)$. Throughout the proof, we shall abuse notation and identify $V_n$ with $P_n A P_n$, $E_n$ with $P_{n+1} A P_n$, and $E_{n}^\top$ with $P_n A P_{n+1}$. Thus, for example, it will be clear that for $x \in S_n$, $E_n \delta_x$ is defined and is a vector in $\ell^2(G)$ with support in $S_{n+1}$. Accordingly, we also regard $\ell^2(S_n)$ as a subspace of $\ell^2(G)$ under the natural inclusion. It is clear that the conclusions of Lemmas \ref{commute} and \ref{commute2} hold for the modified operators as well.

For two vectors $v, w \in \ell^2 (G)$ we write $v \sim w$ when $v= c w$ for some $c \neq 0$.

We shall construct the subspaces $\mathcal{H}_r$ inductively, starting with $\mathcal{H}_0$, which we take to be the cyclic subspace spanned by $A$ and $\delta_o$. First, note that $\delta_o$ is an eigenvector for $\Lambda_{0,+j}$, as well as for  $ (E_{0}^\top E_{1}^\top \ldots  E_{j}^\top V_{j} E_{j} \ldots E_{1} E_0)$, for any $j$. This is simply because the ranges of these operators are contained in $\ell^2(S_0)$ which is one-dimensional. We want to show that the set $\{\phi_0^0, \phi_1^0, \ldots \}$ obtained from Gram-Schmidting $\{\delta_o, A \delta_o, \ldots \}$ has the property described in \emph{(ii)} of the statement of the theorem.

Clearly, $\phi_0^0=\delta_o$. Now, $A \delta_o=V_0 \delta_o+E_0 \delta_o=E_0\delta_o$. If $E_0 \delta_o=0$ then $o$ is disconnected from the rest of $G$ which stands in contradiction to the graph being connected. Thus, $E_0^\top E_0 \delta_o \neq 0$  and it is clear that $\phi_1^0 \sim E_0\delta_o$. We write $E_0^\top E_0 \delta_o=c_0 \delta_o$ for some $c_0 \neq 0$.

Now, $A \phi_1^0 \sim A E_0 \delta_o=E_0^\top E_0 \delta_o+V_1 E_0 \delta_o+E_1 E_0 \delta_o$. We claim that the first two terms in the sum are in the subspace spanned by $\phi_0^0$ and $\phi_1^0$. For the first term this is obvious. For the second term, this follows from
\begin{align*}
V_1 E_0 \delta_o=\frac{1}{c_0}V_1 E_0 (E_0^\top E_0) \delta_o=\frac{1}{c_0}E_0 E_0^\top V_1 E_0 \delta_o
\end{align*}
(which follows from Lemma \ref{commute}) together with the fact that $\delta_o$ is an eigenvector of $E_0^\top V_1 E_0$. Thus, in order to obtain $\phi_2^0, \phi_3^0 ,\ldots$ we need to Gram-Schmidt $\{E_1E_0 \delta_o, A E_1 E_0 \delta_o, \ldots\}$. From this it is clear that if $E_1E_0 \delta_o=0$ we are done. Otherwise, $\phi_2^0 \sim E_1 E_0 \delta_o$ so that  $\supp (\phi_2^0 )\subseteq S_2$. Moreover, there exists $c_1 \neq 0$ such that $\Lambda_{0,+1} \delta_o =c_1\delta_o$.

Now, $A E_1 E_0 \delta_o = E_1^\top E_1 E_0 \delta_o+V_2 E_1 E_0 \delta_o+E_2 E_1E_0\delta_o$. Again, we claim that $\{E_1^\top E_1 E_0 \delta_o, V_2 E_1 E_0 \delta_o\} \subseteq \textrm{span}\{\phi_0^0, \phi_1^0, \phi_2^0\}$. To see that, write
\begin{align*}
V_2 E_1 E_0 \delta_o = \frac{1}{c_1} V_2 E_1 E_0 \Lambda_{0,+1} \delta_o=E_1 E_0 (E_0^\top E_1^\top V_2 E_1 E_0) \delta_o \sim E_1 E_0 \delta_o,
\end{align*}
and
\begin{align*}
E_1^\top E_1 E_0 \delta_o&= E_1^\top E_1 E_0\frac{ (E_0^\top E_0)}{c_0} \delta_o =\frac{1}{{c_0}}\Lambda_{1,+1}\Lambda_{1,-1}E_0 \delta_o =\frac{1}{{c_0}}\Lambda_{1,-1}\Lambda_{1,+1}E_0 \delta_o\\ &=\frac{1}{c_0}(E_0 E_0^\top) (E_1^\top E_1) E_0 \delta_o =\frac{E_0 \Lambda_{0,+1}\delta_o}{c_0} \sim E_0 \delta_o.
\end{align*}

Both equations above follow from the fact that $\delta_o$ is an eigenvector for the appropriate operators, together with the commutation relations from Lemma~\ref{commute}.

Now, suppose we have shown that, for some $n$, $\phi_j^0 \sim E_j E_{j-1} \ldots E_0 \delta_o \neq 0$ for all $j \leq n$, so that there exist nonzero constants $c_0, c_1, c_2 \ldots, c_n$ such that $\Lambda_{0, +j} \delta_o =c_j \delta_o$, $j \leq n$. Consider
\beq \no
A \phi_n^0 \sim E_n^\top E_n \ldots E_0 \delta_o+V_n E_n \ldots E_0 \delta_o+E_{n+1}E_n \ldots E_0 \delta_o,
\eeq
and write
\beq \no
\begin{split}
V_n E_n \ldots E_0 \delta_o &= \frac{1}{c_n}V_n E_n \ldots E_0 \Lambda_{0,+n}\delta_o \\
&=(E_n \ldots E_0) (E_0^\top \ldots E_n^\top V_n E_n \ldots E_0) \delta_o \\
& \sim E_n \ldots E_0 \delta_o
\end{split}
\eeq
and
\beq \no
\begin{split}
 E_n^\top E_n \ldots E_0 \delta_o &=(E_n^\top E_n)E_{n-1}\ldots E_0\frac{\Lambda_{0,+(n-1)}}{c_{n-1}}\delta_o \\
&=\frac{1}{c_{n-1}} \Lambda_{n,+0}\Lambda_{n,-n}E_{n-1}\ldots E_{0}\delta_{o}\\ &=\frac{1}{c_{n-1}} \Lambda_{n,-n}\Lambda_{n,+0}E_{n-1}\ldots E_{0}\delta_{o}\\ &=\frac{1}{c_{n-1}} E_{n-1}\ldots E_0 (E_0^\top\ldots E_{n-1}^\top) E_n^\top E_n (E_{n-1} \ldots E_0)\delta_0 \\
&=\frac{1}{c_{n-1}} E_{n-1}\ldots E_0 \Lambda_{0,+n}\delta_o \\
& \sim  E_{n-1}\ldots E_0 \delta_o,
\end{split}
\eeq
to see that in order to obtain $\{\phi_{n+1}^0, \phi_{n+2}^0,\ldots \}$ we need to Gram-Schmidt $\{E_{n+1}E_n \ldots E_0 \delta_o, A E_{n+1}\ldots E_0 \delta_o\}$. As before, if $E_{n+1}E_{n}\ldots E_0 \delta_o =0$, then we are done. Otherwise, $\phi_{n+1}^0 \sim E_{n+1}E_n \ldots E_0 \delta_o$ and, in particular, $\supp (\phi_{n+1}) \subseteq S_{n+1}$. Moreover, there exists $c_{n+1} \neq 0$ so that $\Lambda_{0, +(n+1)} \delta_o =c_{n+1}\delta_o$. Thus, we have obtained $\mathcal{H}_0$ with the desired properties. Note $n(0)=0$ and $\phi_0^0=\delta_o$.

We now proceed to construct $\mathcal{H}_1$. Let $n(1)$ be the first $n$ so that $\ell^2 (S_n) \nsubseteq \mathcal{H}_0$. Let $W_1$ be the subspace of $\ell^2 (S_{n(1)})$ which is orthogonal to $E_{n(1)-1}\ldots E_0\delta_o$. Then $R_1 \equiv \dim W_1 \geq 1$. First, note that for any $\psi \in W_1$ and any $j < n(1)$, $E_j^\top E_{j+1}^\top \ldots E_{n(1)-1}^\top\psi =0$. This follows from the orthogonality of $W_1$ to $E_{n(1)-1}\ldots E_0\delta_o$, together with the fact that for any $n<n(1)$ $E_n \ldots E_0 \delta_o$ spans $\ell^2(S_n)$.

By Lemmas~\ref{commute} and~\ref{commute2}, we may choose an orthogonal basis $\psi_1, \ldots ,\psi_{R_1}$ of $W_1$, whose members are all mutual eigenvectors of $\Lambda_{n(1), +j}$, $j=1,2,\ldots$, of $V_{n(1)}$ and of $ (E_{n}^\top E_{n+1}^\top \ldots  E_{n+j}^\top V_{n+j} E_{n+j} \ldots E_{n+1} E_n)$, $j=1,2,\ldots$. We now take $\mathcal{H}_1, \mathcal{H}_2,\ldots,\mathcal{H}_{R_1}$ to be the cyclic subspaces spanned by $A$ and $\psi_1, \psi_2,\ldots,\psi_{R_1}$, respectively (so that $n(1)=n(2)=\ldots=n(R_1)$). The proof that these subspaces have the properties required in \emph{(ii)}, is precisely the same as the one for $\mathcal{H}_0$ and uses only the commutation relations together with the fact that the cyclic vectors are eigenvectors for the appropriate operators.

It should be clear at this point how to define $n(R_1+1)$. This will be the first $n$ such that $\ell^2(S_n) \nsubseteq \oplus_{j=0}^{R_1}\mathcal{H}_j$. We let $W_2$ be the subspace of $\ell^2(S_{n(R_1+1)})$ which is orthogonal to  $\oplus_{j=0}^{R_1}\mathcal{H}_j$ and we let $R_2 \geq 1$ be its dimension. Again, we may choose an orthogonal basis to $W_2$ all of whose members are mutual eigenvectors to the appropriate operators and by considering the appropriate cyclic subspaces, construct $\mathcal{H}_{R_1+1}, \ldots \mathcal{H}_{R_1+R_2}$.

Having defined $R_1, R_2, \ldots, R_k$ and $\mathcal{H}_0,\ldots,\mathcal{H}_{\sum_{s=1}^k R_s}$ for some $k \geq 2$, we define $n(R_k+1)$ similarly, to be the distance from the root of the first sphere whose corresponding subspace is not contained in $\mathcal{H}_0 \oplus \ldots \oplus \mathcal{H}_{\sum_{s=1}^k R_s}$. Proceeding as above, we define $R_{k+1}$ and $\mathcal{H}_{\sum_{s=1}^k R_s+1}, \ldots,\mathcal{H}_{\sum_{s=1}^{k+1} R_s}$ as above.

We thus obtain a sequence of subspaces $\mathcal{H}_r$ with the properties described in the theorem. The fact that $\oplus_{r=0}^\infty \mathcal{H}_r \subseteq \ell^2(G)$ follows from the fact that each one of the subspaces is a subspace of $\ell^2(G)$. The reverse inclusion follows since $\{\delta_x \}_{x \in G}$ is a complete orthonormal set in $\ell^2(G)$ and for each $x \in G$, there is some $r_x$ such that $\delta_x \in \oplus_{r=0}^{r_x}\mathcal{H}_r$. Finally, the fact that $A$ decomposes as a direct sum follows immediately from the invariance of each of the subspaces under $A$ and the fact that $A$ is a bounded operator (recall we are assuming that the vertex degrees are bounded).

The proof for $\Delta$ on general strongly path commuting graphs (where the degree is not necessarily bounded), goes along the same lines. The decomposition of $\ell^2(G)$ is precisely the same (where $\widetilde{\ }$ is placed over the appropriate operators). The delicate issue is the fact that $\Delta$ may be unbounded and so we need to consider the issue of selfadjointness for its restrictions, $\Delta_r$, to the various $\mathcal{H}_r$. Moreover, the unboundedness means that $\Delta=\oplus_r \Delta_r$ is also not immediately obvious.

We first claim that, for each $r$, $P_r D(\Delta)=D(\Delta)\cap \mathcal{H}_r$ where $P_r$ is the orthogonal projection onto $\mathcal{H}_r$. The fact that $D(\Delta)\cap \mathcal{H}_r\subseteq P_r D(\Delta)$ is obvious, as well as the fact that $P_r D(\Delta) \subseteq \mathcal{H}_r$. The only nontrivial fact is the fact that $P_r D(\Delta) \subseteq D(\Delta)$. Let $\vphi \in D(\Delta)$ and let $(\vphi^L)$ be a sequence of compactly supported functions satisfying $\vphi^L \rightarrow \vphi$ and $\Delta \vphi^L \rightarrow \Delta \vphi$ as $L \rightarrow \infty$ in $\ell^{2}(G)$. Such a sequence exists since the compactly supported functions are a core for $\Delta$. Clearly $P_r \vphi^L \rightarrow P_r \vphi$. Moreover, by the construction of  $\mathcal{H}_r$, it follows that $\Delta P_r \vphi^L=P_r \Delta \vphi^L$. Thus, $\Delta P_r \vphi^L=P_r \Delta \vphi^L \rightarrow P_r \Delta \vphi$ as $L \rightarrow \infty$, and in particular $\Delta P_r \vphi^L$ has a limit. Since $\Delta$ is closed, it follows that $\lim_{L \rightarrow \infty}P_r \vphi^L=P_r \vphi \in D(\Delta)$. Moreover, this argument also shows that $P_r \Delta \vphi= \Delta P_r \vphi$.

Now, let $\Delta_r=\Delta P_r$ defined on $P_r D(\Delta)$. It immediately follows that $\Delta_r$ is symmetric and closed.
Let $\vphi \in D(\Delta)$ and write $\vphi=\sum_r \vphi_r$ where $\vphi_r=P_r \vphi$. We claim that $\Delta \vphi=\sum_r \Delta_r \vphi_r$. To see this, fix $\varepsilon>0$ and let $L$ be so large that $\|\vphi -\vphi^L\| <\varepsilon$ and also $\|\Delta \vphi-(\Delta \vphi)^L\| < \varepsilon$, where $\psi^L$ is the restriction of $\psi$ to a ball of radius $L$ around the root. By the construction of the $\mathcal{H}_r$ above and since the action of $\Delta$ is local, there exists $R >0$ such that $(\Delta \sum_{r \geq R} \vphi_r)^L=0$. By the invariance of $\mathcal{H}_r$ under $\Delta$, it follows that $\Delta \sum_{r \geq R} \vphi_r$ is orthogonal to $\Delta \sum_{r < R} \vphi_r$ and therefore $\|\Delta \sum_{r \geq R} \vphi_r \| \leq \|\Delta \vphi-(\Delta \vphi)^L\| < \varepsilon$. Thus, $\| \Delta \vphi-\sum_{r<R} \Delta_r \vphi_r\|= \| \Delta  \sum_{r \geq R} \vphi_r \| < \varepsilon$ which implies the claim.

Finally, the preceding paragraph implies that $\Delta_r$ are all selfadjoint, since if $\Delta_r^* \vphi=i \vphi$ then this implies that $\Delta^* \vphi=i \vphi$ which is impossible by the selfadjointness of $\Delta$. Since $\Delta_r$ is symmetric and closed, this ends the proof.
\end{proof}

We now proceed to show that family preserving graphs are strongly path commuting. We need a few lemmas.

\begin{lemma} \label{spherical-symmetry}
Family preserving graphs are spherically symmetric, i.e.,\ if $G$ is a rooted family preserving graph, then for any $n$ and any two vertices $x,y \in S_n$ there is a rooted graph automorphism, $\tau$ such that $\tau(x)=y$.
\end{lemma}

\begin{proof}
We prove this by induction on $n$. It is obviously true for $n=0$, since $S_0=\{o\}$. It is also clearly true for $n=1$ since any $x,y \in S_1$ are backward brothers. Now assume the statement is true for all $j \leq n$, and let $x,y \in S_{n+1}$. If they are backward brothers then we are done. Otherwise, assume that $w \in S_n$ is a neighbor of $x$ and $z \in S_n$ is a neighbor of $y$. Then there is an automorphism, $\tau$ such that $\tau(w)=z$ (by the induction hypothesis). Thus, $\tau(x)$ is a backward brother of $y$ and so there exists $\gamma$ such that $\gamma \tau (x)= y$. We are done.
\end{proof}

\begin{lemma} \label{path-auto}
Let $G$ be a family preserving graph. Let $x, y \in S_n$ and assume there is a $k$-forward path between $x$ and $y$. Then there is a rooted graph automorphism, $\tau$ such that $\tau(x)=y$ and $\tau \upharpoonright S_{n+k}= Id$. Similarly, if  there is a $k$-backward path between $x$ and $y$, then there is a rooted graph automorphism, $\tau$ such that $\tau(x)=y$ and $\tau \upharpoonright S_{n-k}=Id$.
\end{lemma}

\begin{proof}
We prove the forward path case. The backward path case is similar. We use induction on $k$. For $k=1$ this is just the definition of a family preserving graph. Now assume the statement is true for $k$ and suppose there is a $(k+1)$-forward path $(x_0,x_1,\ldots, x_{k+1},\ldots, x_{2k-1}, x_{2k})$ between $x, y \in S_n$ with $x_0=x$ and $x_{2k}=y$.
Then, $(x_1, \ldots, x_{k+1}, \ldots,x_{2k-1})$ is a $k$-forward path between $x_1 \in S_{n+1}$ and $x_{2k-1}\in S_{n+1}$. By the induction hypothesis, there exists an automorphism, $\tau$ such that $\tau(x_1)=x_{2k-1}$ and $\tau \upharpoonright S_{n+k+1} =Id$. It follows that $\tau(x)$ is a forward brother of $y$. Thus, from the property of family preserving graphs, there exists an automorphism, $\gamma$ such that $\gamma(\tau(x))=y$ and $\gamma \upharpoonright S_{n+j}=Id$ for any $j \geq 1$. Then, $\gamma \circ \tau$ is the required automorphism.
\end{proof}

Lemma \ref{path-auto} has two important corollaries:

\begin{corollary} \label{path-correspondence}
Let $G$ be a family preserving graph. Let $x,y \in S_n$ and assume there is a $k$-forward path between $x$ and $y$. Then there is a one-one and onto correspondence between distinct $k$-forward paths between $x$ and $y$ and $k$-forward paths between $x$ and itself.

Similarly, if there is a $k$-backward path between $x$ and $y$ then there is a one-one and onto correspondence between $k$-backward paths between $x$ and $y$ and $k$-bakcward paths between $x$ and itself.
\end{corollary}

\begin{proof}
We prove the $k$-forward path case. By Lemma \ref{path-auto} there is a rooted graph automorphism $\tau$ such that $\tau(x)=y$ and $\tau \upharpoonright S_{n+k}=Id$. Let $(x,x_1,\ldots,x_k,x_{k+1} \ldots, x_{2k-1},x)$ be a $k$-forward path from $x$ to itself. Then $(x, x_1, \ldots,x_k, \tau(x_{k+1}),\ldots,\tau(x_{2k-1}), \tau(x)=y)$ is a $k$-forward path from $x$ to $y$. Clearly, since $\tau$ is an automorphism, this maps the set of paths from $x$ to itself injectively into the set of paths from $x$ to $y$. Using $\tau^{-1}$, we see this map is onto as well.
\end{proof}

\begin{corollary} \label{backward-forward}
Let $G$ be a family preserving graph. Let $x, y \in S_n$. Then: \\
(i) There is a $(k,\ell)$-f.b.\ path from $x$ to $y$ iff there is an $(\ell,k)$-b.f.\ path from $x$ to $y$.\\
(ii) There is a tailed-$k$-forward path from $x$ to $y$ iff there is a headed-$k$-forward path from $x$ to $y$. \\
(iii) There is a tailed-$k$-backward path from $x$ to $y$ iff there is a headed-$k$-backward path from $x$ to $y$.
\end{corollary}

\begin{proof}
To prove $(i)$, we shall show that if there is a $(k,\ell)$-f.b.\ path from $x$ to $y$ then there is an $(\ell,k)$-b.f.\ path from $x$ to $y$. The result will follow by symmetry.

Assume there is a $(k,\ell)$-f.b.\ path from $x$ to $y$. Then there is $z \in S_n$ such that there is a $k$-forward path from $x$ to $z$, and there is an $\ell$-backward path from $z$ to $y$. By Lemma~\ref{path-auto}, there is an automorphism, $\tau$, such that $\tau(z)=x$ and $\tau \upharpoonright S_{n+k}=Id$. Consider $\tau(y)$. Since $\tau$ is an automorphism, there is an $\ell$-backward path from $x=\tau(z)$ to $\tau(y)$ (as there is at least one path from $y$ to a vertex in $S_{n+k}$). But, since $\tau \upharpoonright S_{n+k}=Id$, there is also a $k$-forward path from $\tau(y)$ to $y$. Thus, there is an $(\ell, k)$-b.f.\ path from $x$ to $y$ (passing through $\tau(y)$).

The proofs of $(ii)$ and $(iii)$ are similar.
\end{proof}

\begin{proof}[Proof of Theorem \ref{family-preserving}]
We first want to show that
\beq \no
\#fb_{k,\ell}(x,y)=\#bf_{\ell,k}(x,y)
\eeq
for any $x,y \in S_n$. For any $x \in S_n$, let
\beq \no
F_x^k=\{z \in S_n \mid \textrm{there exists a } k\textrm{-forward path between } x \textrm{ and } z\},
\eeq
and
\beq \no
B_x^\ell=\{z \in S_n \mid \textrm{there exists a }  \ell\textrm{-backward path between } x \textrm{ and } z\}.
\eeq
Let
\beq \no
Z^{k,\ell}_{x,y}=F_x^k \cap B_y^\ell.
\eeq
Then, Corollary~\ref{backward-forward} \emph{(i)} says that $Z^{k,\ell}_{x,y}$ is not empty iff $Z^{k,\ell}_{y,x}$ is not empty. Clearly, if they are both empty then $\#fb_{k,\ell}(x,y)=0=\#bf_{\ell,k}(x,y)$. Thus, we restrict our attention to the case when they are not empty. By Corollary~\ref{path-correspondence}, the number of $k$-forward paths from $x$ to $z \in Z^{k,\ell}_{x,y}$ is the same as the number of $k$-forward paths from $x$ to itself. By the spherical symmetry, this number is independent of $x \in S_n$. Call it $\kappa$. In the same way, the number of $\ell$-backward paths from $y$ to $z \in Z^{k,\ell}_{x,y}$ is a constant. Call it $\lambda$. Clearly,
\beq \no
\#fb_{k,\ell}(x,y)=\kappa \cdot \#Z^{k,\ell}_{x,y} \cdot \lambda
\eeq
and
\beq \no
\#bf_{\ell,k}(x,y)=\kappa \cdot \#Z^{k,\ell}_{y,x} \cdot \lambda.
\eeq
Thus, it suffices to show that $\#Z^{k,\ell}_{x,y}=\#Z^{k,\ell}_{y,x}$.

Now, note that if $z_1, z_2 \in Z^{k,\ell}_{x,y}$ then there is a $k$-forward path from $z_1$ and $z_2$. This is since there is an automorphism $\tau$ such that $\tau (x)=z_2$ and $\tau \upharpoonright S_{n+k} = Id$. Thus, if $(z_1,x_1,\ldots,x_k, x_{k+1},\ldots, x_{2k-1},x)$ is a $k$-forward path between $z_1$ and $x$, then $(z_1,x_1,\ldots, x_k, \tau(x_{k+1}),\ldots, \tau(x)=z_{2})$ is a $k$-forward path between $z_1$ and $z_2$. In the same way, there is an $\ell$-backward path between $z_1$ and $z_2$.

On the other hand, it is easy to show in a similar way that if $z_1 \in Z^{k,\ell}_{x,y}$ and $z_2 \in S_n$ is a vertex such that there exist both a $k$-forward path and an $\ell$-backward path between $z_1$ and $z_2$, then $z_2 \in Z^{k,\ell}_{x,y}$ (i.e., choose $\tau$ such that $\tau(z_{1})=z_{2}$).

Thus, if $z \in Z^{k,\ell}_{x,y}$, then $Z^{k,\ell}_{x,y}$ coincides with the set of vertices in $ S_n$ such that there exist both a $k$-forward path and an $\ell$-backward path between them and $z$. By spherical symmetry, the size of this set is independent of $z \in S_n$. Thus, $\#Z^{k,\ell}_{x,y}=\#Z^{k,\ell}_{y,x}$, and we get
\beq \no
\#fb_{k,\ell}(x,y)=\#bf_{\ell,k}(x,y).
\eeq

To prove
\beq \no
\widetilde{\#f}_k(x,y)=\widehat{\#f}_k(x,y)
\eeq
define $F_k(x,y)$ as the number of vertices in $S_n$ such that there is a $k$-forward path connecting them to $x$ and such that there is an edge connecting them to $y$. Then follow the same strategy as above, to reduce the proof to showing that $\# F_k(x,y)= \# F_k(y,x)$. As before, by Corollary~\ref{backward-forward}, we know they are either both empty or both non-empty, and we assume they are both non-empty.

Let $z \in F_k(x,y)$. Then, as in the proof of Corollary~\ref{backward-forward}, there is an automorphism $\tau$ such that $\tau(z)=x$ and $\tau \upharpoonright S_{n+k} = Id$. Then, $\tau(y)$ is a nearest neighbor of $\tau(z)=x$ and also there is a $k$-forward path between $y$ and $\tau(y)$. Now, consider $\tau^{-1}(F_k(y,x))$. The elements of $F_k(y,x)$ are precisely vertices with a $k$-forward path connecting them to $y$ and an edge connecting them to $x$. Thus, since $\tau \upharpoonright S_{n+k}=Id$ and $\tau^{-1}(x)=z$, the elements of $\tau^{-1}(F_k(y,x))$ are precisely the vertices with a $k$-forward path connecting them to $y$ and an edge connecting them to $z$. In other words $\tau^{-1}(F_k(y,x))=F_k(y,z)$. On the other hand, $F_k(x,y)=F_k(z,y)$ since the elements of $F_k(x,y)$ are precisely the vertices with a $k$-forward path connecting them to $z$ and an edge connecting them to $y$. Thus, we have reduced the problem to showing that $\# F_k(y,z)= \# F_k(z,y)$. The advantage of this is that $y$ and $z$ are neighbors. Thus, by property $(iii)$ of family preserving graphs there is an automorphism $\gamma$ interchanging $z$ and $y$. Clearly, $\gamma$ maps $F_k(y,z)$ onto  $F_k(z,y)$ so the sets have equal size and we are done.

The proof that
\beq \no
\widetilde{\#b}_k(x,y)=\widehat{\#b}_k(x,y)
\eeq
follows the exact same lines of the previous argument. We therefore omit it. It follows that any family preserving graph is path commuting. Moreover, since, by Lemma~\ref{spherical-symmetry}, they are spherically symmetric, we see family preserving graphs are strongly path commuting.
\end{proof}

\section{Antitrees and trees with complete spheres} \label{s:Proof_Applications}
In this section, we provide the proofs for the corollaries for the two classes of examples introduced in Subsection~\ref{ss:antitrees} and~\ref{ss:treeswithcompletespheres}. We start with the decomposition for antitrees and then consider trees with complete spheres.

\subsection{Decomposition for antitrees}

Let an antitree be given. The following theorem shows that the Laplacian on an antitree decomposes into one infinite Jacobi matrix arising from the subspace of spherically symmetric functions and infinitely many one-dimensional Jacobi matrices. In particular, the Laplacian has infinitely many compactly supported eigenfunctions.

\begin{theorem}\label{t:antitree} Let $\Delta$ be the Laplacian on an antitree. Then, $\Delta$ is unitarily equivalent to the operator $$J_{a,b}\oplus\bigoplus_{n\geq1}\bigoplus_{j=1}^{s_{n}-1}M_{n}\quad \mbox{ on } \quad \ell^{2}(\bbN_{0})\oplus\bigoplus_{n\geq1}\bigoplus_{j=1}^{s_{n}-1}\bbC,$$ where $J_{a,b}$ is the Jacobi matrix with off-diagonal $a=(\sqrt{s_{n}s_{n+1}})_{n\geq0}$ and diagonal
$b=(s_{n-1}+s_{n+1})_{n\geq0}$ and $M_{n}$ is the multiplication operator on $\C$ by the number $(s_{n-1}+s_{n+1})$, $n\geq1$. In particular,
\begin{align*}
    \sigma(\Delta) =\sigma(J_{a,b})\cup\{s_{n-1}+s_{n+1}\}_{n\geq1}.
\end{align*}
Moreover, every function that is orthogonal to the spherically symmetric functions and supported on $S_{n}$ is an eigenfunction of $\Delta$ for the eigenvalue  $(s_{n-1}+s_{n+1})$, ${n\geq1}$.
\end{theorem}
\begin{proof} Let $\varphi_{n}$, $n\geq0$, be the functions that take the value $1/\sqrt{s_{n}}$ on  $S_{n}$ and zero elsewhere. Then, $\{\varphi_{n}\}_{n\geq0}$ is an orthonormal basis for the subspace of compactly supported spherically symmetric functions. Clearly, $\Delta$ maps a  spherically symmetric function to a spherically symmetric function. Calculating
\begin{align*}
\widehat{a}_{n}&=\langle{\Delta \varphi_{n+1},\varphi_{n}}\rangle=\frac{1}{\sqrt{s_{n}}}\sum_{x\in S_{n}}\Delta \varphi_{n+1}(x)=\frac{-1}{\sqrt{s_{n}s_{n+1}}}\sum_{x\in S_{n}} \sum_{y\in S_{n+1}}1 =-\sqrt{s_{n}s_{n+1}},\\
b_{n}&=\langle{\Delta\varphi_{n},\varphi_{n}}\rangle =\frac{1}{\sqrt{s_{n}}}\sum_{x\in S_{n}}\Delta\varphi_{n}(x) =\frac{1}{s_{n}}\sum_{x\in S_{n}}(s_{n-1}+s_{n+1}) ={s_{n-1}+s_{n+1}},
\end{align*}
yields that $\Delta$ restricted to the spherically symmetric functions in the domain of $\Delta $ is unitarily equivalent to the Jacobi matrix $J_{\widehat{a},{b}}$ on $\ell^{2}(\bbN_{0})$. Moreover, $J_{\widehat{a},b}$ is unitarily equivalent to $J_{-\widehat{a},b}=J_{a,b}$.

It can be directly checked that any function that is orthogonal to the spherically symmetric functions and supported on $S_{n}$ is an eigenfunction with the eigenvalue $(s_{n-1}+s_{n+1})$, $n\geq1$.
Since any function orthogonal to the spherically symmetric function can be decomposed into functions that are supported on the spheres, we are done.
\end{proof}

\textbf{Remark.}
From the proof of the theorem above we can derive similar statements for the normalized Laplacian and the adjacency matrix (in the case where it defines a self adjoint operator). Let an antitree be given and $(s_{n})$ be the cardinalities of the spheres.

(a) The normalized Laplacian $\widetilde{\Delta}$ on $\ell^{2}(G,\deg)$ is the bounded operator given by $(\widetilde{\Delta}\varphi)(x)=\deg(x)^{-1}\sum_{y\sim x}(\varphi(x)-\varphi(y))$. Then, $\sigma(\widetilde{\Delta})=\sigma(J_{\widetilde a},1)\cup\{1\}$, where the Jacobi matrix $J_{\widetilde a,1}$ has constant diagonal $1$ and $\widetilde a_{n}=a_{n}/\sqrt{b_{n}b_{n+1}}=\sqrt{s_{n}s_{n+1}/((s_{n-1}+s_{n+1})(s_{n}+s_{n+2}))}$. Moreover, $1$ is an eigenvalue of  infinite multiplicity.

(b) Suppose the restriction of the adjacency matrix $A$ to the finitely supported functions is an essentially selfadjoint operator. Then, $\sigma(A)=\sigma(J_{a,0})\cup \{0\}$  where the eigenvalue zero has infinite multiplicity and $J_{a,0}$ is the Jacobi matrix with zero diagonal and $a$ is  as in the theorem.
\medskip

The statement of  Theorem~\ref{c:antitrees1}
follows directly from  Theorem~\ref{t:antitree} above.  Theorem~\ref{c:antitrees2} is a consequence of Remling's theorem.

\begin{proof}[Proof of Theorem~\ref{c:antitrees2}] By \cite[Theorem~1.1]{Rem} the Jacobi matrix $J_{a,b}$ has absolutely continuous spectrum if and only if the sequences $a$ and $b$ are eventually periodic. This is the case if and only if $(s_{n})$ is eventually periodic: Let $q$ be the period of $a$ and $b$ for $n\geq N$. If $s_{n}=s_{n+qk}$ for some $n\geq N,k\ge1$, then it follows by induction using $a_{j}=a_{j+qk}$, $j\geq N$ that  $(s_{n})$ is eventually periodic with period $qk$. Since $a$, $b$ are eventually periodic,  $(s_{n})$ takes only finitely many values. Therefore, there exists  $n$ and $k$ such that  $s_{n}=s_{n+qk}$. The other direction is trivial.
\end{proof}

\subsection{Decomposition for trees with complete spheres}
Let us turn to the decomposition for trees with complete spheres. Recall that such a graph $G(k,\gamma)$ is a tree $T(k)$ with branching sequence $k\in\N^{N}$ and within the sphere $S_{n}$ all vertices are connected by edges whenever $\gamma_{n}=1$, $n\ge1$.

\begin{theorem}\label{t:tree1} Let $k\in \N^{\N}$,  $\gamma\in\{0,1\}^{\N}$ and the Laplacian  $\Delta$ on $G(k,\gamma)$ be given.
Then, $\Delta$ is unitarily equivalent to the direct sum of operators $$\bigoplus_{l\geq0}\bigoplus_{j=1}^{m_{l}}J_{a^{(l)},b^{(l)}}\quad \mbox{ on } \quad \bigoplus_{l\geq0}\bigoplus_{j=1}^{m_{l}}\ell^{2}(\bbN_{0}),$$
where $J_{a^{(0)},b^{(0)}}$ is the Jacobi matrix with diagonal with $b^{(0)}_{0}=k_{1}$, $b^{(0)}_{n}=k_{n+1}+1$, $n\geq 1$,
and off-diagonal $a^{(0)}_{n}=\sqrt{k_{n+1}}$, $n\geq 0$ and $m_{0}=1$. Furthermore,
$J_{a^{(l)},b^{(l)}}=J_{(a^{(0)}_{n+l}),(b^{(0)}_{n+l})}+v^{(l)}$, where $v^{(l)}$ is the diagonal matrix with entries
$v^{(l)}_{n}=v_{n+l}$ with $v_{0}=0$, $v_{n}=\gamma_{n}\prod_{j=1}^{n}k_{j}$, $n\geq1$, and $m_{n}=\prod_{j=1}^{n}k_{j}-\prod_{j=1}^{n-1}k_{j}$, $l\geq 1$, ${n\geq0}$. In particular,
\begin{align*}
    \sigma(\Delta)=\bigcup_{l\geq 0}\sigma(J_{a^{(l)},b^{(l)}}).
\end{align*}
\end{theorem}
\begin{proof}We obtain the decomposition into invariant subspaces exactly in the same way as it was done  for trees in \cite{AF,Br,GG}: For $l=0$ we take the orbit of the delta function at the root and, for $l\geq 0$, we take the orbits of $m_{l}$ orthogonal functions supported on $S_{l}$ whose values on vertices with common father sum up to zero. We proceed by 'Gram-Schmidting' the functions in each orbit to get an orthonormal basis of functions that are supported on one sphere each. For $l\geq 1$, all these functions are zero on $S_{j}$, $j<l$, and the functions are constant
on  $S_{j}\cap T_{x}$, $j>l$, where $T_{x}$ is the forward tree  of some $x\in S_{l}$. Moreover, summing up over any sphere gives zero for any such function. Thus, a direct computation  as in the proof of Theorem~\ref{t:antitree} gives the particular values of $a^{(l)}$, $b^{(l)}$.
\end{proof}

Let us turn to the proof of the corollaries for these graphs.
\begin{proof}[Proof of  Theorem~\ref{c:treeswithcompletespheres1}] The Jacobi matrix associated to the subspace of spherically symmetric functions is eventually periodic and thus has finitely many bands of absolutely continuous spectrum and  finitely many discrete eigenvalues outside of the bands. For any other Jacobi matrix $v^{(l)}_{n}=\prod_{j=l}^{n+l}k_{j}\to\infty$, $n\to\infty$,  so, the spectrum is purely discrete. Since the bottom of the spectrum of $J_{a^{(l)},b^{(l)}}$ is larger or equal than $v^{(l)}_{1}$, there only finitely many eigenvalues of $\Delta$ smaller than $v^{(l)}_{1}$. Since $v^{(l)}_{1}=v_{l}\to\infty$,  $l\to\infty$, the spectrum of $\Delta $ on the orthogonal complements of the spherically symmetric functions is purely discrete.
\end{proof}

\begin{proof}[Proof of  Theorem~\ref{c:treeswithcompletespheres2}] Since $\Delta$ is a positive operator, there is no spectrum below zero.
 Moreover,  it is clear that $0$ is not an eigenvalue: suppose there is a function $\phi\in\ell^{2}(G)$ with $\Delta\phi=0$. Since $0=\langle{\phi,\Delta\phi\rangle}=\frac{1}{2}\sum_{x \sim y}|\phi(x)-\phi(y)|^{2}$ and $\phi\in\ell^{2}(G)$ it follows $\phi\equiv 0$.

The Jacobi matrix $J_{a^{(0)},b^{(0)}}$ associated to the subspace of spherically symmetric functions is the free Laplacian on $\N$ with a finite rank perturbation. Hence, its absolutely continuous spectrum is $[0,4]$ and there are at most finitely many eigenvalues outside $[0,4)$.

The Jacobi matrices  $J_{a^{(l)},b^{(l)}}+v^{(l)}$ associated with the  subspaces orthogonal to the spherically symmetric functions have no absolutely continuous spectrum by \cite{Rem} since the potentials $v^{(l)}_n=v_{n+ l}$ with $v_n=\gamma_{n}(\kappa-1)$, $l=1,\ldots,\kappa-1$, are not eventually periodic and $v^{(l)}$ take finitely many values. Moreover, there is no point spectrum in $(0,4)$ by  the arguments in the proof of Theorem~4.1  in \cite{SiSt}. The essential spectrum of these matrices is $[0,4] \cup \{(2+\sqrt{4+\kappa^2})\}$ by \cite[Theorem 1.7]{LS1} and a simple computation. It follows that the singular continuous spectrum fills the interval $[0,4]$ and that ${4}$ might be an eigenvalue and an accumulation point of eigenvalues, and that $(2+\sqrt{4+\kappa^2})$ is an accumulation point of eigenvalues (and perhaps an eigenvalue as well).
\end{proof}

\end{document}